\magnification 1250
\mathcode`A="7041 \mathcode`B="7042 \mathcode`C="7043
\mathcode`D="7044 \mathcode`E="7045 \mathcode`F="7046
\mathcode`G="7047 \mathcode`H="7048 \mathcode`I="7049
\mathcode`J="704A \mathcode`K="704B \mathcode`L="704C
\mathcode`M="704D \mathcode`N="704E \mathcode`O="704F
\mathcode`P="7050 \mathcode`Q="7051 \mathcode`R="7052
\mathcode`S="7053 \mathcode`T="7054 \mathcode`U="7055
\mathcode`V="7056 \mathcode`W="7057 \mathcode`X="7058
\mathcode`Y="7059 \mathcode`Z="705A
\def\spacedmath#1{\def\packedmath##1${\bgroup\mathsurround
=0pt##1\egroup$} \mathsurround#1
\everymath={\packedmath}\everydisplay={\mathsurround=0pt}} 
\def\nospacedmath{\mathsurround=0pt
\everymath={}\everydisplay={} } \spacedmath{2pt}
\def\qfl#1{\buildrel {#1}\over \longrightarrow}
\def\hfl#1#2{\normalbaselines{\baselineskip=0pt
\lineskip=10truept\lineskiplimit=1truept}\nospacedmath\smash
{\mathop{\hbox to 12truemm{\rightarrowfill}}
\limits^{\scriptstyle#1}_{\scriptstyle#2}}}
\def\diagram#1{\def\normalbaselines{\baselineskip=0truept
\lineskip=10truept\lineskiplimit=1truept}   \matrix{#1}}
\def\vfl#1#2{\llap{$\scriptstyle#1$}\left\downarrow\vbox
to 6truemm{}\right.\rlap{$\scriptstyle#2$}}
\font\eightrm=cmr8      \font\sixrm=cmr6
\catcode`\@=11
\def\eightpoint{%
 \textfont0=\eightrm \scriptfont0=\sixrm
\scriptscriptfont0=\fiverm
  \def\rm{\fam\z@\eightrm}%
   \rm}\catcode`\@=12

\def\note#1#2{\footnote{\parindent
0.4cm$^#1$}{\vtop{\eightpoint\baselineskip12pt\hsize15.5truecm
\noindent #2}}\parindent 0cm}
\def\iso{\vbox{\hbox to .8cm{\hfill{$\scriptstyle\sim$}\hfill}
\nointerlineskip\hbox to .8cm{{\hfill$\longrightarrow $\hfill}} }}

\def\pc#1{\tenrm#1\sevenrm}
\def\up#1{\raise 1ex\hbox{\smallf@nt#1}}
\def\tx{\kern-1.5pt -}
\def\cqfd{\kern 2truemm\unskip\penalty 500\vrule height
4pt depth 0pt width 4pt}  \def\virg{\raise .4ex\hbox{,}}  
\def\decale#1{\smallbreak\hskip 28pt\llap{#1}\kern 5pt}
\def\no{n\up{o}\kern 2pt}
\def\ind{\par\hskip 0,8truecm\relax}

\parindent=0cm
\def\moins{\mathrel{\hbox{\vrule height 3pt depth -2pt
width 6pt}}}
\def\rond{\kern 1pt{\scriptstyle\circ}\kern 1pt}

\def\rk{\mathop{\rm rk}\nolimits}
\def\Hom{\mathop{\rm Hom}\nolimits}

\def\Pic{\mathop{\rm Pic}\nolimits}
\frenchspacing
\font\cc=cmcsc10
\input amssym.def
\input amssym
\vsize = 25truecm
\hsize = 16truecm
\font\tensan=cmssdc10
\font\sevensan=cmssdc10 at 7pt
\font\fivesan=cmssdc10 at 5pt
\newfam\sanfam
 \textfont\sanfam=\tensan \scriptfont\sanfam=\sevensan
  \scriptscriptfont\sanfam=\fivesan
  \def\san{\fam\sanfam\tensan}%
\def\ext{{\san \Lambda}}
\voffset = -.5truecm
\baselineskip15pt
\overfullrule=0pt
\null
\centerline{\bf Some stable vector bundles  with reducible theta
divisors}\smallskip
\centerline{Arnaud {\cc Beauville}}
\vskip1truecm
{\bf Introduction}\smallskip 
\ind Let
$C$ be a curve of genus $g$, and 
$E$  a vector bundle on
$C$, of rank $r$;  assume that the {\it slope}  $\mu:={1
\over r}\deg E$ of
$E$ is an integer. Let $J^{\nu}$ be the
variety  (isomorphic to the Jacobian of $C$) parametrizing line
bundles of degree $\nu:=g-1-\mu$ on $C$. We  say that
$E$ {\it admits a theta divisor} if 
$H^0(E\otimes L)=0$  for  $L$ general in $J^\nu$.
 If this is the case, the  locus
$$\Theta_{E}=\{ L\in J^{\nu}\ |\ H^0(E\otimes
L)\not=0\}$$has a natural structure of effective
divisor in $J^{\nu}$,  the
{\it theta divisor} of $E$. Its class in $H^2(J^\nu,{\bf Z})$ is
$r\theta$, where $\theta\in H^2(J^\nu,{\bf Z})$ is the
class of the principal polarization.
 This (generalized) theta divisor plays a key
 role in the recent work on vector bundles on curves --
see for instance  [B] for an overview. 
\ind If $E$ admits a theta divisor, it is semi-stable (otherwise $E$
contains a sub-bundle $F$ of slope $>\mu$, and by
Riemann-Roch $H^0(F\otimes L)$, and therefore $H^0(E\otimes
L)$, is non-zero for all $L\in J^{\nu}$). The converse does
not hold, at least in rank $\ge 4$: Raynaud has constructed
examples of stable vector bundles with no theta divisor [R].
Further examples have been constructed recently by Popa [P].
\ind If $E$ is semi-stable but not stable, its theta divisor (if it
exists) is not integral: more precisely, $E$ admits a filtration with
stable quotients $E_1,\ldots,E_p$, and we have $\Theta
_E=\Theta _{E_1}+\ldots +\Theta _{E_p}$. One may ask,
conversely, if the reducibility of $\Theta _E$ implies  that $E$
is not stable. A counter-example has been given by Raynaud
(unpublished), who constructed a rank 2 stable vector bundle on
a curve of genus 3 with reducible theta divisor. Such an
example can only occur on a special curve, because in rank 2 the
divisor $\Theta _E$ characterizes the vector bundle
$E$ [B-V], and on a general curve the only reducible divisors on
$J^\nu$ with cohomology class $2\theta$ are the theta divisors
of rank 2 decomposable bundles.

\ind We describe in this note a counter-example of a different
nature, namely a family of stable  vector bundles of rank $g$
which exist on any  curve of genus $g$. They are defined
by the exact sequence
$$0\rightarrow M_L \longrightarrow H^0(C,L)\otimes_{\bf
C}{\cal O}_C\qfl{\rm ev_L} L\rightarrow 0$$where $L$ is a line
bundle generated by its global sections, and ${\rm ev_L}$ the
evaluation map. These  vector bundles have been intensively
studied, notably by Green and Lazarsfeld (see in particular [L]),
Paranjape and Ramanan [P-R], and more recently in [P] and
[F-M-P]. In the latter paper the authors determine the theta
divisor of $M_K$ and of its exterior powers; we will do the
same here in the case   of a line bundle $L$ of degree
$2g$ (so that $M_L$ has rank $g$). We will prove, in a
somewhat more precise form, the following result:
\smallskip 
{\bf Theorem}$.-$  {\it Let $C$ be a non-hyperelliptic curve, and
$L$ a sufficiently general line bundle of degree $2g$ on $C$. The
vector bundle $M_L$ and its exterior powers
$\ext^2M_L,\ldots,\ext^{g-1}M_L $  are stable and
admit a reducible theta divisor.}
 \smallskip 
\ind  An interesting extra feature of our examples is that
there exists a semi-stable, decomposable vector bundle on $C$
with the same theta divisor as $M_L$; thus in rank $\ge 3$ the 
divisor $\Theta_E$ does not characterize the bundle $E$ any
more.
\vskip1truecm
{\bf Notation}
\ind We fix a curve $C$ of genus $g$; except in Remark 2
below, we assume throughout that $C$ {\it is not hyperelliptic}.
For $d\in{\bf Z}$, we denote by $J^d$ the translate of the
Jacobian of $C$ parametrizing line bundles of degree $d$ on
$C$, and by
$C_d$ the locus of effective divisor classes in $J^d$. If $p,q\in
{\bf Z}$ the difference variety $C_p-C_q$ lies in $J^{p-q}$.
\vskip1truecm
{\bf I. The Theta divisor of ${\bf E_L}$}\smallskip 
\ind Let $L$ be a line bundle of degree $2g$ on the curve $C$.
It is spanned by its global sections, so we have an 
exact sequence
$$0\rightarrow M_L\longrightarrow H^0(L)\otimes_{\bf C}{\cal
O}_C\longrightarrow L\rightarrow 0\ ,$$where $M_L$ is a rank
$g$ vector bundle. We put $E_L:=M_L^*$. 
\ind Though this will not be used  in the sequel, let us recall the
geometric interpretation of $E_L$. Let
$\varphi$ be the morphism of $C$ into the
projective space  ${\bf P}:={\bf
P}(H^0(L))$ defined
by the linear system $|L|$; in view of the Euler exact sequence
$$0\rightarrow {\cal O}_{\bf P}\longrightarrow
H^0(L)^*\otimes_{\bf C}{\cal O}_{\bf P}(1) \longrightarrow T_{\bf
P}\rightarrow 0\ ,$$we have $E_L=\varphi^*T_{\bf
P}\otimes L^{-1} $.  
\ind The vector bundle $E_L$ has rank
$g$ and determinant $L $, hence slope $2$.

\smallskip 
{\bf Proposition 1}$.-$ a) {\it The vector bundle 
$E_L$ has a theta divisor} 
$$\Theta _{E_L}=(C_{g-2}-C)+\Theta
_{L\otimes K^{-1}}\quad{\it in}\ J^{g-3}\ .$$
\ind b) {\it $E_L$ is semi-stable; it is stable  if and only
if $L$ is very ample.}
\smallskip 
{\it Proof} : We will first compute set-theoretically the theta
divisor $\Theta _{M_L}$ of $M_L$. By definition this is  the set of
line bundles $P\in J^{g+1}$ such that the multiplication map
$\mu:H^0(L)\otimes H^0(P)\rightarrow H^0(L\otimes P)$ 
is not injective. Let us distinguish three cases:
\ind (i)  If $h^0(P)>2$ we have $\dim
H^0(L)\otimes H^0(P)>\dim H^0(L\otimes P)$, thus $P\in
\Theta_{M_L}$. 
\ind (ii) Assume that $h^0(P)=2$ and that the pencil $|P|$ has a 
base point.   Both spaces $H^0(L)\otimes H^0(P)$
and $H^0(L\otimes P)$ have the same dimension $2g+2$. If
$\mu$ is injective,  it is surjective, and
the linear  system
$|L\otimes P|$ has a base point; this  is impossible since
$\deg (L\otimes P)= 3g+1$. Thus we have again 
$P\in\Theta_{M_L}$.
  
\ind (iii) Finally assume that $|P|$ is a base-point free pencil. From
the exact sequence
$$0\rightarrow P^{-1} \longrightarrow
H^0(P)\otimes_{\bf C} {\cal O}_C\longrightarrow P\rightarrow 0$$
we get an exact sequence
$$0\rightarrow H^0(L\otimes P^{-1}) \longrightarrow
H^0(L)\otimes_{\bf C} H^0(P)\qfl{\mu}
H^0(L\otimes P)\ ;$$thus $\mu$ is not injective in that case if
and only if $H^0(L\otimes P^{-1})\not= 0$.
\ind The line bundles $P$ in case (i) and (ii) are exactly those
which can be written
$P'(x)$, for some point $x$ of $C$ and some line bundle $P'$ in
$J^g$ with
$h^0(P')\ge 2$; the ones in case (iii) are those of the form
$L\otimes P'^{-1}$, with $P'\in\Theta \i J^{g-1}$.
Since
$\Theta _{E_L}$ is the image of $\Theta _{M_L}$ by the
isomorphism
$L\mapsto K\otimes L^{-1} $ of $J^{g+1}$ onto $J^{g-3}$, we
obtain (still set-theoretically) $\Theta _{E_L}=(C_{g-2}-C)\cup
\Theta _{L\otimes K^{-1} }$. Now $C_{g-2}-C$ is an irreducible
divisor with cohomology class $(g-1)\theta$ (see e.g. [F-M-P],
Prop. 3.7), and $\Theta _{L\otimes K^{-1} }$ is a (ordinary) theta
divisor; since
$\Theta _{E_L}$ has cohomology class $g\theta$, we get the
equality a).
\ind  Since $E_L$ admits a theta divisor, it is semi-stable. 
Moreover, its stable components are
$L':=L\otimes K^{-1} $ and a rank $(g-1)$ bundle. If $E_L$ is not
stable,
$L'$ is either a sub- or a quotient bundle of $E_L$. The latter
case cannot occur since
$E_L$ is generated by its global sections and $L'$ is not. 
 Now using the exact sequence
$$0\rightarrow L^{-1}\otimes L'^{-1} \longrightarrow
H^0(L)^*\otimes_{\bf C}  L'^{-1} \longrightarrow E_L\otimes
L'^{-1} \rightarrow 0$$and Serre duality we see that
$\Hom(L',E_L)$ is zero if and only if the multiplication map
$H^0(L)\otimes H^0(L)\rightarrow H^0(L^2)$ is surjective, that
is, $L$ is normally generated [G-L]. By [G-L], Thm. 1, this is the
case if and only if $L$ is very ample.\cqfd
\medskip
{\it Remarks}$.-$ 1) If $L$ is not very ample, we have
$L=K(D)$, with $D$ an effective divisor of degree $2$. The
snake lemma applied to the commutative diagram
$$\diagram{0 \longrightarrow & M_K  &\longrightarrow
& H^0(K)\otimes {\cal O}_C &\longrightarrow &K
&\longrightarrow 0\cr &\vfl{}{} && \vfl{}{} && \vfl{}{}\cr
0 \longrightarrow & M_L &\longrightarrow
& H^0(L)\otimes {\cal O}_C &\longrightarrow
&L&\longrightarrow 0 }$$provides an exact sequence
$0\rightarrow M_K\rightarrow M_L\rightarrow {\cal
O}_C(-D)\rightarrow 0$; thus $E_L$ is an extension of
$E_K$ by ${\cal O}_C(D)$. This extension is
non-trivial, as we already observed that ${\cal
O}_C(D)$ cannot be a quotient of $E_L$.\smallskip 
\ind 2) If $C$ is hyperelliptic, the divisor $C_{g-2}-C$ is equal to
$\Theta_H $, where $H$ is the hyperelliptic pencil on $C$. By
specialization we get $\Theta _{E_L}=(g-1)\Theta_H
+\Theta_{L\otimes K^{-1} }$. The line bundle $L$ is not linearly
normal [L-M], so $E_L$ is not stable.
\medskip
\ind The difference variety $C_{g-2}-C$ is the theta divisor of the
bundle $E_K$ [P-R]; therefore:\par
 {\bf Corollary 1}$.-$ {\it Assume that $L$ is very ample. The
stable bundle $E_L$ and the decomposable bundle $E_K\oplus
(L\otimes K^{-1}) $ have the same theta
divisor}.\cqfd\smallskip 
\ind The equality still holds of course when $L$ is not very
ample, but becomes immediate, since in that case the second
bundle is the sum of the stable components of the first one.
\ind In view of the results of [B-N-R], this corollary can be
rephrased as follows. Let ${\cal SU}_C(g)$ be the moduli space of
semi-stable rank $g$ vector bundles on $C$ with trivial
determinant, and let  ${\cal L}$ be the positive generator of
$\Pic ({\cal SU}_C(g))$ (the {\it determinant bundle}). Let
$B_{\cal L}$ be the base locus of the linear system $|{\cal L}|$.
\smallskip {\bf Corollary 2}$.-$ {\it The map $\varphi _{\cal L}:
{\cal SU}_C(g)\moins B_{\cal L}\,\longrightarrow\, {\bf
P}(H^0({\cal L}))$ defined by the line bundle ${\cal L}$ is not
injective}.
\ind Indeed this map can be identified with the map which
associates to a vector bundle its theta divisor [B-N-R]. Twisting
 $E_L$ and $E_K\oplus (L\otimes K^{-1})$ by  a line bundle
$\lambda$ on $C$ with $\lambda^{- g}= L$, we get two
different points of ${\cal SU}_C(g)\moins B_{\cal L}$
with the same image under $\varphi _{\cal L}$.\cqfd

\vskip1truecm

{\bf II. The Theta divisor of ${\bf \ext}^{p}
{\bf E_L}$}\smallskip 
\ind We now consider the exterior power $\ext^pE_L$; this is a
vector bundle of rank ${g\choose p}$ and slope $2p$, so its
theta divisor, if it exists, lie in $J^{g-1-2p}$.
\smallskip 
{\bf Proposition 2}$.-$
{\it Let $1\le p\le g-1$. If $L$ is general enough,
 the vector bundle $\ext^pE_L$ is stable and
admits a theta divisor}\note{1}{The second term is the
translate of $C_{g-p}-C_{p-1}\i J^{g+1-2p}$ by the element
$K\otimes L^{-1}$ of $J^{-2}$.}
$$\Theta _{\ext^pE_L}= (C_{g-p-1}-C_{p})+
(C_{g-p}-C_{p-1}+K\otimes L^{-1} )\ .$$
\smallskip 
{\it Proof} : We first prove that $\ext^pE_L$ admits a theta
divisor when $L$ is general enough. Since this is an open
property, it is sufficient to prove this for a particular choice of
$L$: we take $L=K(D)$, with $D$ an effective divisor of degree
$2$. The exact sequence
$$0\rightarrow {\cal O}_C(D)\rightarrow E_L\rightarrow
E_K\rightarrow 0$$ (Remark 1) gives rise to an exact
sequence
$$0\rightarrow \ext^{p-1}E_K\,(D)\rightarrow
\ext^pE_L\rightarrow \ext^pE_K\rightarrow 0\ .$$Since by
[F-M-P] each exterior power $\ext^qE_K$ admits a theta divisor,
so does
$\ext^pE_L$.
\ind Let us now prove that the theta divisor $\Theta _{E_L}$,
when it exists, is given by the formula of the Proposition. The
divisor
$C_q-C_{g-1-q}$ has cohomology class
${g-1\choose q}\theta$ ([F-M-P], Prop. 3.7), so both sides
of the formula  have  cohomology
class ${g\choose p}\theta$. It suffices therefore to prove that
each component of the right hand side is contained in $\Theta
_{\ext^pE_L}$.

\ind As in [P] and [F-M-P], we will use the following observation
of Lazarsfeld [L]: if $x_1,\ldots,x_{g-1}$ are generic points of
$C$, there is an exact sequence
$$0\rightarrow \bigoplus_{i=1}^{g-1}{\cal
O}_C(x_i)\longrightarrow E_L\longrightarrow L(-\sum
x_i)\rightarrow 0
$$which gives rise as above to an exact sequence
$$0\rightarrow \bigoplus_{i_1<\ldots<i_p} {\cal
O}_C(x_{i_1}+\ldots+x_{i_p})\longrightarrow
\ext^pE_L\longrightarrow \bigoplus_{j_1<\ldots<j_{g-p}}
L(-x_{j_1}-\ldots-x_{j_{g-p}})\rightarrow 0\ .$$
This gives:
\ind $\bullet$ $H^0(\ext^pE_L(-x_1-\ldots-x_p))\not= 0$, hence
the inclusion $C_{g-p-1}-C_{p}\i \Theta _{\ext^pE_L}$;
\ind $\bullet$ $H^0(\ext^pM_L\otimes L
(-x_1-\ldots-x_{g-p}))\not= 0$, hence $H^0(\ext^pM_L\otimes L
(-D))\not= 0$ for all $D$ in $C_{g-p}-C_{p-1}$; by Serre duality
this gives $H^0(\ext^pE_L\otimes K\otimes L^{-1} (D))\not= 0$, 
hence the inclusion $C_{g-p}-C_{p-1}+K\otimes L^{-1} \i \Theta
_{\ext^pE_L}$.
\ind It remains to prove that $\ext^pE_L$ is stable. Since $L$ is
generic, $E_L$ is stable (Proposition 1), so $\ext^pE_L$ is {\it
polystable} -- that is, direct sum of stable bundles with the same
slope $2p$. If $\ext^pE_L$ is not stable for $L$ generic, it is
decomposable for all values of $L$; we will see that this is not the
case when  $L$ is of the form
$K(D)$, with $D$ effective of degree $2$. In that case we have  by
Remark 1 an exact sequence
$$0\rightarrow \ext^{p-1}E_K\,(D)\longrightarrow
\ext^pE_L\longrightarrow \ext^pE_K\rightarrow 0$$
where $ \ext^{p-1}E_K\,(D)$ and $\ext^pE_K$ are stable with
slope $2p$; if $\ext^pE_L$ is decomposable, this exact sequence
splits. The following easy lemma shows that this is not the case,
and thus concludes the proof of the Proposition.\cqfd
\smallskip 
{\it Lemma}$.-$ {\it Let $X$ be a scheme over a field of
characteristic $0$, and let 
 $$0\rightarrow M
\rightarrow E\rightarrow F\rightarrow 0\leqno({\cal E})$$ be a
non-split exact sequence of vector bundles on $X$, with $\rk
M=1$. The associated exact sequences
$$0\rightarrow \ext^{p-1}F\otimes M\longrightarrow
\ext^pE\longrightarrow \ext^pF\rightarrow 0\leqno(\ext^p{\cal
E})$$ do not split for $1\le p\le \rk F$}.
\smallskip {\it Proof} : Let $i:F^*\otimes M\rightarrow
{\cal H}om(\ext^pF,\ext^{p-1}F\otimes M)$ be the linear map
deduced from the interior product. A straightforward
computation shows  that the  class  of the
extension $(\ext^p{\cal E})$ in 
$H^1(X,{\cal H}om(\ext^pF,\ext^{p-1}F\otimes M))$  is the image 
by $H^1(i)$ of  the class of the extension  $({\cal E})$ in
$H^1(X,F^*\otimes M)$. But in characteristic zero $i$ admits a
retraction $c^{-1} \rho$, where $c={\rk F-1\choose p-1}$ and
$\rho:{\cal H}om(\ext^pF,\ext^{p-1}F\otimes M)\rightarrow
F^*\otimes M$ is the map deduced from the interior product
$\ext^pF^*\otimes
\ext^{p-1}F\rightarrow F^*$. Thus $H^1(i)$ is injective, and the
lemma follows.\cqfd
\medskip

 \ind As in section I this gives:\par
{\bf Corollary 1}$.-$ {\it The vector bundles $\ext^pE_L$ and
$\ext^pE_K\oplus (\ext^{p-1}E_K\otimes L\otimes K^{-1}) $
have the same theta divisor. In particular, if $L$ is general
enough, the map
$\varphi _{\cal L}: $\break ${\cal SU}_C({g\choose p})\moins
B_{\cal L}\,\longrightarrow\, {\bf P}(H^0({\cal L}))$ defined by
the line bundle
${\cal L}$ is not injective.}
\medskip
\ind Let us conclude by a link with the main theme of
[F-M-P], the so-called {\it minimal resolution conjecture} for
the curve $C$ embedded into ${\bf
P}^g:={\bf P}(H^0(L))$. We have to refer to [F-M-P] for the 
statement of the conjecture, which is a bit technical. Let us just
say that it describes, for all general finite subsets $\Gamma \i C$
of cardinality $\ge g+1$,  the minimal graded resolution of
the ideal $I_\Gamma $ of $\Gamma$ in the coordinate ring
$S={\bf C}[X_0,\ldots,X_g]$ of
${\bf P}^g$. By
 Corollary 1.8 of [F-M-P], this conjecture holds if and only if
each of the bundles $\ext^pE_L$ admit  a theta
divisor. Thus:\smallskip 
 {\bf Corollary 2}$.-$ {\it The curve $C$,  embedded into ${\bf
P}^g$ by a general linear system of degree $2g$, satisfies the
``minimal resolution conjecture" in the sense of} [F-M-P].\cqfd
\vskip2cm
\centerline{ REFERENCES} \vglue15pt\baselineskip12.8pt
\def\num#1#2#3{\smallskip\item{\hbox to\parindent{\enskip
[#1]\hfill}}{\cc #2}: {\sl #3}.}
\parindent=1.3cm 
\num{B} {A. Beauville} {Vector bundles on curves and
generalized theta functions: recent results and open problems}
Current topics in complex algebraic geometry, 17--33, Math. Sci.
Res. Inst. Publ. {\bf 28}, Cambridge Univ. Press
(1995).
\num{B-N-R}{A. Beauville, M.S. Narasimhan, S. Ramanan}
{Spectral curves  and the generalised theta divisor} J. Reine
Angew. Math. {\bf 398} (1989), 169--179.
\num{B-V}{S. Brivio, A. Verra} {The theta divisor
of ${\rm SU}\sb C(2,2d)\sp s$ is very ample if $C$ is not
hyperelliptic} Duke Math. J. {\bf 82} (1996),  503--552. 
\num{F-M-P} {\ G. Farkas, M. Musta\c{t}\v{a}, M. Popa}{Divisors 
on ${\cal M}_{g,g+1}$ and the Mini\-mal Resolution Conjecture
for points on canonical curves} Preprint {\tt
math.AG/0104187}, to appear in Duke Math. J.
\num{G-L}{M. Green, R. Lazarsfeld}{On the projective
normality of complete  linear series on an algebraic curve}
Invent. Math. {\bf 83} (1986), 73--90.
\num{L}{R. Lazarsfeld}{A sampling of vector bundle techniques
in the study of linear series}  Lectures on Riemann Surfaces,
World Scientific Press, Singapore (1989), pp. 500--559.
\num{L-M}{H. Lange, G. Martens}  {Normal generation and
presentation of line bundles of low degree on curves} J. Reine
Angew. Math. {\bf 356} (1985), 1--18.
\num{P}{M. Popa}{On the base locus of the generalized
theta divisor}   C. R. Acad. Sci. Paris {\bf 329} (1999), S\'erie
I, 507--512.

\num{P-R}{K. Paranjape, S. Ramanan}{On the
canonical ring of a curve} Algebraic geometry and commutative
algebra, vol. II, 503--516, Kinokuniya, Tokyo (1988).

\num{R}{M. Raynaud} {Sections des fibr\'es vectoriels sur
une courbe} Bull. Soc. Math. France {\bf 110} (1982),
103--125.

\vskip1cm
\def\pc#1{\eightrm#1\sixrm}
\hfill\vtop{\eightrm\hbox to 5cm{\hfill Arnaud {\pc BEAUVILLE}
\hfill}
 \hbox to 5cm{\hfill Institut Universitaire de France\hfill}\vskip-2pt
\hbox to 5cm{\hfill \&\hfill}\vskip-2pt
 \hbox to 5cm{\hfill Laboratoire J.-A. Dieudonn\'e\hfill}
 \hbox to 5cm{\sixrm\hfill UMR 6621 du CNRS\hfill}
\hbox to 5cm{\hfill {\pc UNIVERSIT\'E DE}  {\pc NICE}\hfill}
\hbox to 5cm{\hfill  Parc Valrose\hfill}
\hbox to 5cm{\hfill F-06108 {\pc NICE} Cedex 02\hfill}}
\end